\documentclass[12pt]{amsart}
\usepackage{amsmath, amssymb, graphicx}
\usepackage{amssymb}
\usepackage{amsthm}
\usepackage[arrow,matrix]{xy}
\newtheorem{Thm}{Theorem}
\newtheorem{Lem}[Thm]{Lemma}

\newtheorem{Prop}[Thm]{Proposition}
\newtheorem{Rm}{Remark}
\newtheorem{Def}{Definition}

\def\C{{\mathbb C}}

\def\E{{\mathbb E}}
\def\R{{\mathbb R}}

\def\P{{\mathbb P}}
\def\J{{\mathbb J}}

\def\Z{{\mathbb Z}}
\def\W{{\mathbb W}}

\def\O{{\mathbb O}}
\def\H{{\mathbb H}}
\def\V{{\mathbb V}}

\def\Z{{\mathbb Z}}

\def\CG{\mathcal G}

\def\CA{\mathcal A}

\def\CL{\mathcal L}

\def\s2x{\hbox{$S^2 \times S^2$}}

\def \Di {D\!\!\!\!/\,}

    \def\sqr#1#2{{\vcenter{\hrule height.#2pt
            \hbox{\vrule width.#2pt height#1pt \kern#1pt
            \vrule width.#2pt}\hrule height.#2pt}}}
    \def\square{\mathchoice\sqr67\sqr67\sqr{2.1}6\sqr{1.5}6}

\def\qed{~\hfill$\square$}

\begin{document}

\title[]{Deformations in $G_2$ Manifolds }
\author{Selman Akbulut and Sema Salur}
\thanks{First named author is partially supported by NSF grant DMS 0505638}
\keywords{mirror duality, calibration}
\address{Department  of Mathematics, Michigan State University, East Lansing, MI, 48824}
\email{akbulut@math.msu.edu }
\address {Department of Mathematics,  University of Rochester, Rochester, NY, 14627 }
\email{salur@math.rochester.edu }
\date{\today}

\begin{abstract}
Here we  study the deformations of associative submanifolds inside  a $G_2$ manifold $M^7$ with a calibration $3$-form $\varphi $. A choice of $2$-plane field $\Lambda$ on $M$ (which always exits) splits the tangent bundle of $M$ as a direct sum of a $3$-dimensional associate bundle and a complex $4$-plane bundle $TM={\bf E}\oplus {\bf V}$, and this helps us to relate the deformations to Seiberg-Witten type equations.  Here all the surveyed results as well as the new ones about $G_2$ manifolds are proved  by using only the cross product operation (equivalently  $\varphi$).  We feel that mixing various different local  identifications of the rich $G_2$ geometry (e.g. cross product, representation theory and the algebra of octonions) makes  the study of $G_2$ manifolds looks harder then it is (e.g. the proof of McLean's theorem \cite{m}). We believe the approach here makes things easier and keeps the  presentation  elementary. This paper is essentially self contained.

\end{abstract}

\maketitle

\setcounter{section}{0}
\vspace{-0.3in}

\section{${\mathbf G_2}$ manifolds}

 We first review the basic results about $G_2$ manifolds, along the way we give a self contained proof of the McLean's theorem and its generalization \cite{m}, \cite{as1}. A $G_2$ manifold $(M,\varphi, \Lambda )$ with an oriented $2$-plane field  gives various  complex structures on some of  subbundles of $T(M)$. This imposes  interesting structures on the deformation theory  of  its associative submanifolds.  By using this we relate them  to the Seiberg-Witten type equations. 

\vspace{.05in}

Let us  recall some basic definitions (c.f. \cite{b1}, \cite{b2},\cite{hl}):  Octonions give an 8 dimensional division algebra $\O=\H\oplus l \H= \R^8$ generated by  $\langle 1, i, j, k, l, li ,lj, lk\rangle $. The imaginary octonions $im \O =\R^7$ is equipped with the cross
product operation $\times:  \R^7\times \R^7 \to \R^7$ defined  by
$u\times v=im(\bar{v} .u)$. The exceptional Lie group $G_{2}$ is
the linear automorphisms of  $im \O$ preserving this cross product.  It can also be defined in terms of the orthogonal
$3$-frames:
\begin{equation}
G_{2} =\{ (u_{1},u_{2},u_{3})\in (im \O)^3 \;|\; \langle u_{i},u_{j} \rangle=\delta_{ij}, \; \langle
u_{1} \times u_{2},u_{3} \rangle =0 \;  \}.
\end{equation}

\vspace{.1in}

\noindent Alternatively,  $G_2$ is the
subgroup of $GL(7,\R)$ which fixes a particular $3$-form
$\varphi_{0} \in \Omega^{3}(\R ^{7})$, \cite{b1}. Denote
$e^{ijk}=dx^{i}\wedge dx^{j} \wedge dx^{k}\in \Omega^{3}(\R^7)$,
then
$$ G_{2}=\{ A \in GL(7,\R) \; | \; A^{*} \varphi_{0} =\varphi_{0}\; \}. $$

{\vspace{-.15in}

\begin{equation}
\varphi _{0}
=e^{123}+e^{145}+e^{167}+e^{246}-e^{257}-e^{347}-e^{356}.
\end{equation}

{\Def A smooth $7$-manifold $M^7$ has a {\it $G_{2}$ structure} if
its tangent frame bundle reduces to a $G_{2}$ bundle.
Equivalently, $M^7$ has a {\it $G_{2}$ structure} if there is  a
3-form $\varphi \in \Omega^{3}(M)$  such that  at each $x\in  M$
the pair $ (T_{x}(M), \varphi (x) )$ is  isomorphic to $(T_{0}(
\R^{7}), \varphi_{0})$ (pointwise condition). We call
$(M,\varphi)$ a manifold with a $G_2$ structure.}

\vspace{.1in}

A $G_{2}$ structure $\varphi$ on $M^7$ gives an orientation $\mu_{\varphi }=\mu
\in \Omega^{7}(M)$ on $M$,   and $\mu$ determines  a metric
$g=g_{\varphi }= \langle \;,\;\rangle$ on $M$, and a cross product operation  $TM\times TM\mapsto TM$: $(u,v)\mapsto u\times v=u\times_{\varphi} v $  defined as follows: Let
$i_{v}=v\lrcorner $ be the interior product with a vector $v$,
then
\begin{equation*}
\langle u,v \rangle=[ (u\lrcorner\; \varphi ) \wedge (v \lrcorner\; \varphi )\wedge
\varphi  ]/6\mu .
\end{equation*}
\begin{equation}
\varphi (u,v,w) = (v\lrcorner\; u\lrcorner \;\varphi) (w)=\langle u\times v,w \rangle .
\end{equation}



{\Def A manifold with $G_{2}$ structure $(M,\varphi)$  is called a
{\it $G_{2}$ manifold} if   the holonomy group of the Levi-Civita
connection (of the metric $g_{\varphi }$) lies inside of  $G_2 $. In this case $\varphi$ is called integrable.
Equivalently  $(M,\varphi)$ is a $G_{2}$ manifold if $\varphi $ is
parallel with respect to the metric $g_{\varphi }$, that is
$\nabla_{g_{\varphi }}(\varphi)=0$; which is in  turn equivalent to
$d\varphi=0 $, $\;d(*_{g_{\varphi}}\varphi)=0$ (i.e. $\varphi$
harmonic). Also equivalently,  at each point $x_{0}\in M$ there is
a chart  $(U,x_{0}) \to (\R^{7},0)$ on which $\varphi $ equals to
$\varphi_{0}$ up to second order term, i.e. on the image of $U$,
$\varphi (x)=\varphi_{0} + O(|x|^2)$}.

\vspace{.05in}

{\Rm One important class of $G_2$ manifolds are the ones obtained
from Calabi-Yau manifolds. Let $(X,\omega, \Omega)$ be a complex
3-dimensional Calabi-Yau manifold with K\"{a}hler form $\omega$
and a nowhere vanishing holomorphic 3-form $\Omega$, then
$X^6\times S^1$ has holonomy group $SU(3)\subset G_2$, hence is a
$G_2$ manifold. In this case $ \varphi$= Re $\Omega + \omega
\wedge dt$. Similarly, $X^6\times \mathbb{R}$ gives a noncompact
$G_2$ manifold. }

\vspace{.05in}

{\Def Let $(M, \varphi )$ be a $G_2$ manifold. A 4-dimensional
submanifold $X\subset M$ is called {\em coassociative } if
$\varphi|_X=0$. A 3-dimensional submanifold $Y\subset M$ is called
{\em associative} if $\varphi|_Y\equiv vol(Y)$; this condition is
equivalent to $\chi|_Y\equiv 0$,  where $\chi \in \Omega^{3}(M,
TM)$ is the  tangent bundle valued 3-form given by:}
\begin{equation}
\langle \chi (u,v,w) , z \rangle=*\varphi  (u,v,w,z).
\end{equation}

\noindent Equivalence of these  conditions follows from  the `associator
equality' of  \cite{hl}
\begin{equation*}
\varphi  (u,v,w)^2 + |\chi (u,v,w)|^2/4= |u\wedge v\wedge w|^2.
\end{equation*}

\noindent Sometimes $\chi$ is  also called the triple cross product operation and denoted by $\chi(u,v,w) =u\times v\times w$. By imitating the definition of $\chi$,  we can view the usual cross product  as a tangent bundle 2-form  $\psi \in \Omega^{2}(M, TM)$ defined by
\begin{equation}
\langle \psi (u,v) , w \rangle=\varphi  (u,v,w).
\end{equation}

\noindent  As  in the case of  $\varphi$, $\chi $ can  be expressed in terms of in cross product and metric
\begin{equation}
\chi(u,v,w)= -u\times (v\times w)-\langle u,v\rangle w +\langle u,w\rangle v
\end{equation}

\vspace{.1in}

\noindent (c.f. \cite{h}, \cite{hl}, \cite{k}). From (6)  and the identity $u\times v=u.v +\langle u,v\rangle $, the reader can easily check that $2\chi(u,v,w)=(u.v).w- u.(v.w) $, which shows that the associative submanifolds of 
$(M,\varphi)$ are the manifolds where the octonion multiplication of the tangent vectors is ``associative''.  

\vspace{.1in}

We call a $3$-plane $E\subset TM$ {\it associative plane} if  $\varphi|_{E}=vol (E) $, so associate submanifolds $Y^3$ are  submanifolds whose tangent planes are associative. From (2) and (3) we see that  an associative $3$-plane  $E\subset TM$ is a plane generated by three orthonormal vectors in the form $\langle u,v,u\times v \rangle$; and also  if  $V=E^{\perp}$ is its orthogonal complement (coassociative), the cross product induces maps: 
\begin{equation}
 E\times  V \to V,   \mbox{ and} \;\; V\times V\to  E , \mbox{ and} \;\; E\times E \to  E.
 \end{equation}


Note that (4) implies that the $3$-form $\chi$ assigns a normal vector to every oriented $3$-plane in $T(M)$, \cite{as1}, which is zero on the associative planes. Therefore, we can view  $\chi$ as a section of  the $4$-plane bundle $\V=\E^{\perp} \to G_{3}(M)$ over the  Grassmannian bundle of orientable $3$-planes in $T(M)$, where 
$\V$ is orthogonal bundle to the  canonical bundle $\E\to G_{3}(M)$. 
In particular, $\chi$ gives a normal vector field on all oriented  $3$-dimensional submanifolds $f:Y^3 \hookrightarrow(M,\varphi)$, which is zero if the submanifold is associative. This gives an interesting first order flow $\partial f/\partial t =\chi(f_{*}vol ( Y) )$ (which is called $\chi$-flow in \cite{as2}), which appears to push 
$f(Y)$ towards associative submanifolds. 

\vspace{.1 in}

Finally, a useful fact which will be used later is the following: The  $SO(3)$-bundle ${\E}$ is the reduction of the $SO(4)$-bundle $\V$ by the projection 
to the first factor  $ SO(4)=(SU(2)\times SU(2))/{\bf Z}_{2}\to SU(2)/{\bf Z}_2 = SO(3)$, i.e. $\E=\Lambda^{2}_{+}\V$.

\vspace{.1 in}

\section{2-frame fields of ${\mathbf G_2}$ manifolds}

\vspace{.1 in}

 By a theorem of Emery Thomas,  all orientable $7$-manifolds admit non-vanishing $2$-frame fields \cite{t}, in particular they admit non-vanishing oriented $2$-plane fields. Using this, we get a useful additional structure on the tangent bundle of $G_2$ manifolds.
 
 {\Lem A non-vanishing oriented $2$-plane field $\Lambda $ on a manifold with $G_2$-structure $(M,\varphi)$ induces a splitting  of $T(M)={\bf E}\oplus {\bf V}$, where ${\bf E}$ is a bundle of associative $3$-planes, and ${\bf V}={\bf E}^{\perp}$ is a bundle of coassociative $4$-planes. The unit sections $\xi$  of the bundle ${\bf E}\to M$  give complex structures  $J_{\xi} $ on ${\bf V}$.  
 
 \proof Let $\Lambda = \langle u,v \rangle $ be the $2$-plane spanned by the  basis vectors of  an orthonormal  $2$-frame  $\{u,v\}$  in $M$. Then we define ${\bf E}=\langle u,v,u\times v  \rangle$, and  ${\bf V}=  {\bf E}^{\perp}$. We can define the complex structure  on 
 ${\bf V}$ by $J_{\xi} (x)=x\times \xi $.   \qed

\vspace{.1 in}
Similar complex structures   were studied in \cite{hl}. The complex structure  $J_{\Lambda} (z)=\chi (u,v,z)$  of \cite{as1}  turns out to coincide with $J_{v\times u}$ because by (6):
\begin{equation}
\chi(u,v,z)=\chi(z,u,v)=
-z\times (u\times v)-\langle z, u\rangle v + \langle z, v\rangle u=J_{v\times u}(z).
\end{equation}

 $J_{\xi}$ also defines a complex structure on the bigger bundle $\xi^{\perp}\subset TM$.  So it  is natural to study  manifolds   $(M,\varphi, \Lambda )$, with a $G_2$ structure $\varphi$,   and  a nonvanishing oriented $2$-plane field $\Lambda $  inducing the splitting $T(M)= {\bf E}\oplus {\bf V}$, and  ${\bf J}=J_{v\times u}$. Note that each of these terms depend on  $\varphi$ and $\Lambda $.

\vspace{.05 in}

{\Def We call $Y^3\subset (M,\Lambda )$ 
a {\it $\Lambda$-spin} submanifold if $\Lambda |_{Y}\subset TY$, and  call 
$Y^3\subset (M,\varphi, \Lambda )$  a $\Lambda$-associative  submanifold if 
${\bf E}|_{Y}=TY $. }

\vspace{.05 in}

Clearly $\Lambda$-associative submanifolds $Y\subset (M,\varphi, \Lambda )$ are $\Lambda$-spin. Also since  $Y$ has a natural metric induced from the metric of $(M,\varphi)$, we can identify the set of $Spin^c$ structures $Spin^{c}(Y) \cong H^2(Y,\Z)$ on $Y$ by  the homotopy classes of $2$-plane fields on $Y$ (as well as the homotopy classes of vector fields on $Y$). So, any  $\Lambda$-spin submanifold $Y$  inherits a   natural $Spin^c$  structure $s=s(\Lambda)$ from $\Lambda$.  How abundant are the $\Lambda$-associative (or $\Lambda$-spin) submanifolds? Some answers:
 
 \vspace{.05 in}

{\Lem Let $M^7$ be an orientable $7$-manifold, then  
every $Spin^c$ submanifold $(Y^3,s) \subset M^7$  is  $\Lambda$-spin for some $\Lambda$ with $s=s(\Lambda)$, and every  associative $Y\subset (M, \varphi )$   is $\Lambda$-associative for some $\Lambda $. }

\proof Let  $s=\langle u',v'\rangle$ be the $Spin^c$ structure generated by   an orthonormal frame field on $TY$.   By using \cite{t} we choose a nonvanishing  orthonormal $2$-frame field $\{ u,v \}$ on $M$.  Let $V_{2}({\bf R}^7)\to V_{2}(M)\to M$ be the Steifel bundle of $2$-frames in $T(M)$. Now the restriction of this bundle to $Y$ has two sections $\{u',v'\}$ and  $\{u,v\}|_{Y}$ which  are homotopic, since the fiber  $V_{2}({\bf R}^7)$ is $4$-connected.  By the homotopy extension property $\{u',v'\}$ extends to orthonormal $2$-frame field $\{u'',v''\}$ to $M$,  then  we let $\Lambda=\langle u'',v'' \rangle$. Furthermore when $Y$ is associative, we can start with an orthonormal $3$-frame  of $TY$ of the form $ \{ u',v', u' \times v'  \}$, then get the corresponding ${\bf E}_{\Lambda} = \langle u'',v'', u''\times v'' \rangle $, which makes $Y$  $\Lambda$-associative. 
 \qed
 
  \vspace{.1 in}

  More generally, for any manifold with a $G_2$ structure $(M, \varphi)$  we can study the bundle of oriented $2$-planes $G_{2}(M)\to M$ on $M$, and construct the corresponding universal  bundles  ${\E}\to G_2(M)$ and ${\V}\to G_2(M)$, and a complex structure ${\J }$  on ${\V}$, where $\J=J_{\Lambda}$ on the fiber over  $\Lambda =\langle u,v\rangle$. Then each $(M,\varphi, \Lambda )$ is a section of  $G_{2}(M)\to M$, inducing  ${\bf E}, {\bf V}, {\bf J}$. 
We can do the same construction on the bundle of oriented $2$-frames $V_{2}(M)\to M$ and get the same quantities,  in this case we get a {\it hyper-complex structure} on $\V$, i.e. we get three complex structures ${\J}={\J}_1,  {\J}_2,{\J}_3$ on $\V$ corresponding to $J_{u\times v}, J_{u}, J_{v}, $ over each fiber $\{u,v \}$, and  they anti-commute and cyclically commute e.g. ${\J}_1{\J}_{2}={\J}_3 $. Notice also that ${\J}_1$ depends only on the oriented $2$-plane field, whereas   
${\J}_2,{\J}_3$  depend on the $2$-frame field.

 \vspace{.1 in}

\noindent  By using $\J_1$ (or one of the other $\J_p$, $p=2,3$) we can split $\V_{\C}=\W\oplus \bar{\W}$, as a pair of conjugate $\C^2$-bundles ($\pm i$ eigenspaces of $\J_1$). This gives a complex line bundle $K=\Lambda ^{2}(\bar{W})$ which correponds to the $2$-plane field $\Lambda$. Corresponding to $K$ we get a canonical $Spin^c$ structure on $\V$. More specifically, recall that
$ U(2)=(S^1\times S^3) /\Z_{2}$, $SO(4)=(S^3\times S^3)/\Z_{2}$, $Spin^{c}(4)=
(S^3\times S^3\times S^1)/\Z _{2}$, 
\begin{equation}
\begin{array}{ccc}
  &   &   Spin^{c}(4)\\
  &  \nearrow &  \downarrow \\
 U(2) &  \to &   SO(4)\times S^1 
\end{array}
\end{equation}

 \noindent where the horizontal map $[\lambda,A]\mapsto ([\lambda,A], \lambda^2)$  canonically lifts to the map $[\lambda,A]\mapsto (\lambda,A, \lambda)$, where the transition functions $\lambda^2$ corresponds to $K$ (see for example \cite{a}). 
 This means that there are pair of complex $\C^2$-bundles, ${\W}^{\pm}\to V_{2}(M)$ with $\V_{\C}={\W}^{+}\otimes {\W}^{-}$. This fact can be checked directly by taking  ${\W}^{+}=K^{-1}+\C $ and  ${\W}^{-}=\bar{W} $ 
 (note $ \Lambda^{2}(\W) \otimes  \bar{W} \cong \Lambda^{2}(\W)\otimes {\W}^{*}   \cong \W $). This gives an action $\E=\Lambda^{2}_{+}(\V) : \W^{+}\to \W^{+}$; in our case this action will come from cross product structure, Lemma 3 will do this by identifying $\W^{+}$ with $S$.

 \vspace{.1 in}

Note also that from (6) an  (7) the cross product operation $\rho(a)(w)= a\times w$ induces a Clifford representation by  $ \rho  (u\times v) =-\J_{1}$, $\rho (u ) =-\J_{2}$, $\rho (v) =-\J_{3}$
\begin{equation}
\rho :\E \to End(\V).
\end{equation}

\vspace{.1 in}

\subsection {$G_2$ frame fields on $G_2$ manifolds}
$\;$
 
 \vspace{.1 in}

In the case of a manifold with $G_2$ structure $(M,\varphi )$, Thomas's theorem can be strengthen to the conclusion that $M$ admits  a $2$-frame field $\Lambda $, with the property that on the tubular neighborhood of the $3$-skeleton of $M$,  $\Lambda$ is the restriction of a $G_2$ frame field.  To see this, we start with an orthonormal  $2$-frame field $\{u_1,u_2\}$, then let $\Lambda =\langle u_1,u_2, u_1\times u_2 \rangle$ and $\V\to M$ be the corresponding universal $4$-plane bundle as in the last section. Then we pick a unit section $u_3$ of   $\V\to M$ over the $3$-skeleton $M^{(3)}$; there is no obstruction doing this since we are sectioning an $S^3$-bundle over the $3$-skeleton of $M$. Now, from the definition of $G_2$ in (1) we see that $\{u_1,u_2,u_3\}$ is a $G_2$ frame on  $M^{(3)}$.  

 \vspace{.04 in}

{\Def  We call $(M,\varphi, \Lambda)$ a {\it framed $G_2$ manifold}  if $\Lambda$ is the restriction of a $G_2$ frame field on $M$.}

 \vspace{.05 in}

The above discussion says that  every $(M, \varphi)$ admits a $2$-frame field $\Lambda$ such that $(M^{(3)},\varphi, \Lambda)$ is a $G_2$-framed manifold.
 From now on, the notation  $(M,\varphi, \Lambda)$ will refer to a manifold with a $G_2$ structure and a $2$-frame field $\Lambda$, such that on $M^{(3)}$, $\Lambda$  is the restriction of a $G_2$ frame as above.  From the above discussion,  the last condition is equivalent to picking a   nonvanishing section  of $\V \to M^{(3)}$ (called $u_3$ above). This will be useful when studying local deformations of associative submanifolds $Y^3\subset M $ (they live near $M^{(3)}$). Using the same notations of the last section we state:

 \vspace{.05 in}

{\Lem Let $(M,\varphi, \Lambda)$ be a framed $G_2$ manifold. Then we can decompose $\V_{\C}=S\oplus \bar{S}$ as a pair of bundles, each of which is isomorphic to  $W^{+}=K^{-1}+\C$,  and  the cross product $\rho$ induces a representation $\rho_{\C} :\E_{\C} \to End(S) $ given by: 

 \begin{equation*} u\times v\mapsto 
\left(
\begin{array}{ccc}
- i  &   0   \\
 0 &  i    \\     
\end{array}
\right) \;\;\;
u\mapsto 
\left(
\begin{array}{ccc}
 0 & 1     \\
- 1 &   0   \\     
\end{array}
\right) \;\;\;
v\mapsto 
\left(
\begin{array}{ccc}
0 &  i   \\
  i & 0     \\     
\end{array}
\right) 
\end{equation*}
}
\proof We choose a local orthonormal frame $\{e_1,..,e_7\}$ which $\varphi $ is in the form (2) with $\{u\times v,u,v \}=\{e_1,e_2,e_3\}$ (because of the canonical metric we will not distinguish the notations of local frames and coframes). From (2) and (3) we compute the cross product operation, $\J_1,\J_2,\J_3$, and $\W$ from the tables below

\begin{table}[h]
\begin{center}
\begin{tabular}{|c||c|c|c|c|} \hline
$ \times $ & $e_4$  &  $e_5$ &  $e_6 $ &  $e_7$\\
\hline \hline
$e_1$ & $e_5$ & $-e_4$ & $e_7$ & $-e_6$ \\
$e_2$ & $e_6$  & $-e_7$ & $-e_4$ & $e_5$ \\
$e_3$ &  $-e_7$ & $-e_6$ & $e_5$ &  $e_4$\\
\hline
\end{tabular}
\end{center}
\end{table}

\vspace{-.2in}
 \[ \J_{1}:
\left(
\begin{array}{l}
 e_4 \mapsto -e_5       \\
 e_5\mapsto e_4       \\
  e_6\mapsto -e_7     \\
 e_7\mapsto e_6       
\end{array}
\right),
\;\;\;
\J_{2}:
\left(
\begin{array}{l}
 e_4 \mapsto -e_6     \\
 e_5\mapsto e_7       \\
  e_6\mapsto e_4    \\
 e_7\mapsto -e_5      
\end{array}
\right),
\;\;\;
\J_{3}:
\left(
\begin{array}{l}
 e_4 \mapsto e_7     \\
 e_5\mapsto e_6       \\
  e_6\mapsto -e_5    \\
 e_7\mapsto - e_4       
\end{array}
\right)
\]

\vspace{.1 in}

 $\W=\langle e_p -i\; \J_{1}(e_p)\;|\: p=4,..,7 \rangle =\langle e_4 +ie_5, e_6 +ie_7 \rangle_{\C} =\langle E_1, E_2\rangle _{\C}. $

\vspace{.05 in}

$\bar{\W}=\langle e_p +i\; \J_{1}(e_p)\;|\: p=4,..,7 \rangle =\langle e_4 -ie_5, e_6 -ie_7 \rangle_{\C} =\langle \bar{E}_1, \bar{E}_2\rangle _{\C}. $

 \vspace{.1 in}

  $\J_2 \mbox{ and } \J_3  :\W\to \bar{\W}$ are given by $(E_1, E_2)\mapsto (-\bar{E_2},\bar{E_1})$ and $ (i\bar{E_2},-i\bar{E_1})$ respectively;  by composing them with complex conjugation we can view them as complex structures on  $\W $  (hence we get a quaternionic structure on $\W$). We can  decompose $\V_{\C}= S \oplus \bar{S} $, where  $S=\langle {E}_1, J_{2} E_1\rangle _{\C}= \langle {E}_1, -\bar{E}_2\rangle _{\C} $, and hence $\bar{S}= \langle {E}_2, J_{2} E_2\rangle _{\C}=
  \langle E_2, \bar{E_1} \rangle _{\C} $, then it is straightforward to check that,  the maps $\J_p$ give complex structures on $S$  and $\rho (e_p)$ are given by the matrices in the statement of this Lemma, for $p=1,2,3$.  
  
 \vspace{.05 in}
   
  Since $a(e_4+ie_5)+b(e_6-ie_7)= [a(e_4+ie_5)\wedge (e_6+ie_7)]+ b]\otimes (e_6-ie_7)$,  we can identify $S\cong K^{-1} +\C$, i.e.  tensoring with the section $s:=(e_6-ie_7)$ gives the isomorphism. Here s is a nonvanishing section of $\V_{\C}$ which is determined by the unit section $u_{3}$ coming from the $G_2$ framing (discussed above). This is because we can choose $\{ e_4 =u_3 , e_5 =J_{1} (u_3),  e_6 =J_{2} (u_3),  e_7 =J_{3} (u_3) \}$. \qed
   
  
 \vspace{.2 in}
 
There is also the useful bundle map $\sigma:  S \to \E$ induced by  
\begin{equation}
\sigma (z,w) =( \frac{ |z|^2-|w|^2 }{2} ) u\times v+ Re (z\bar{w}) u + Im (w\bar{z}) v.
\end{equation}
 This is the quadratic map which appears in Seiberg-Witten theory, after identifying $\E$ with the Lie algebra $su(2)$ (skew adjoint endomorphisms of $\C^2$ with the inner product given by the Killing form) we get 
  \[ \sigma(x) =\sigma(z,w)= 
\left(
\begin{array}{ccc}
 \frac{ |z|^2-|w|^2 }{2}&    z\bar{w}  \\
  w\bar{z}&    \frac{ |w|^2-|z|^2 }{2}  \\
  \end{array}
\right).
\]

\begin{equation}
\langle \sigma (x), x \rangle = 2 |\sigma(x)|^2=\frac{1}{2} |x|^4 .
\end{equation} 
These identifications are standard tools used Seiberg-Witten theory (c.f \cite{a}).

  \vspace{.2 in}

\subsection{Deforming $G_2$ structures}

$\;$
\vspace{.1 in}

 For a $7$-manifold with a $G_2$ structure $(M,\varphi)$, the space of all $G_{2}$ structures on  $M$ is identified with an open subset of $3$-forms $\Omega^{3}_{+}(M)\subset \Omega^{3}(M)$, which is the orbit of $\varphi$ by the gauge transformations of  $T(M)$.  The orbit is open by the dimension reason (recall that the action of $GL(7,{\bf R})$ on $\Omega^{3}(X)$ has $G_2$ as the stabilizer).  The structure of $\Omega^{3}_{+}(M)$ is nicely explained in \cite{b2} as follows:  By definition, $\Omega^{3}_{+}(M)$ is the space of sections of a bundle over $M$ with fiber $GL(7,{\bf R})/G_2$ (which is homotopy equivalent to $\R\P^7$). Furthermore, the subspace of the $G_2$ structures  inducing the same metric can be parametrized with the space of sections of the  bundle $\R\P^7\to P(T^{*}M\oplus \R)\to M$ with fibers $SO(7)/G_2 =\R\P^7$, where 
 $ P(T^{*}M\oplus \R)$ is the projectivization of $T^{*}(M)\oplus \R$. That is, if $\lambda=[a,\alpha ] $ with $a^2+\alpha^2=1$,  then the corresponding 
 $\varphi_{\lambda}\in\Omega^{3}_{+}(M)$ is  
 \begin{equation}\varphi_{\lambda}=\varphi -2 \alpha^{\#}\lrcorner \;[\;a(*\varphi)+\alpha\wedge \varphi \;]\end{equation}
 
 \noindent  where $\alpha^{\#}$ is the metric dual of $\alpha$.
This is  given in \cite{b2}, written slightly differently. Therefore,  if we start with an integrable $G_2$ structure with harmonic $\varphi $, the space of integrable $G_2$ structures inducing the same metric are parametrized by the sections $\lambda=[a,\alpha ] $, such that  $d\theta=d(*\theta) =0$,  where $ \theta=\alpha^{\#}\lrcorner \;[\;a(*\varphi)+\alpha\wedge \varphi \;]$ and    $* \theta=\alpha \wedge [\;a \varphi - (\alpha^{\#} \lrcorner * \varphi ) \; ]$. It is a natural question whether a submanifold $Y^3\subset (M,\varphi )$ is associative. The following says that any $Y$ can be made associative in 
$(M, \varphi_{\lambda})$, after deforming $\varphi $ to $\varphi_{\lambda}$.

{\Prop Let $(M^7, \varphi )$  be a manifold with a $G_2$ structure, then any  $Spin^c $ submanifold $(Y^3 ,s) \subset M^7$  is a  $\Lambda $-associative submanifold of $(M, \varphi_{\lambda}, \Lambda)$ for some choice of 
$\lambda=[a,\alpha]$ and a plane field $\Lambda $.}

\proof By Lemma 2, we can assume $Y$ is $\Lambda$-spin for some $\Lambda =\langle u,v\rangle $. Hence this gives an orthogonal  splitting $T(M)={\bf E}\oplus {\bf V}$, with ${\bf E}=\langle u,v, u\times v\rangle $.  Choose a unit vector field $w$ in $Y$ orthogonal to  $\langle u,v\rangle |_{Y}$, then extend $w$ to $M$. Now we want to choose $\lambda=[a,\alpha]$ so that if $(u\times v)_{\lambda}$ is the cross product corresponding to the $G_2$ structure $\varphi_{\lambda}$, then 
$(u\times v)_{\lambda}|_{Y}=w$. 

\vspace{.05in}

By (13), and  the rules  
$(u\times v)^{\#} =v \lrcorner\; u \; \lrcorner \;\varphi $ and 
$(u\times v)^{\#} _{\lambda}=v \lrcorner\; u \; \lrcorner \;\varphi_{\lambda} $  we get 
\begin{eqnarray*}  (u\times v)_{\lambda}/2&=& (u\times v)/2 -|\alpha |^2 (u\times v)  - a \chi(u,v, \alpha^{\#})  \\ && +\alpha(v)(u\times \alpha^{\#} ) -\alpha(u)(v\times \alpha^{\#} )+\alpha^{\#} \langle u\times v,\alpha^{\#} \rangle .
\end{eqnarray*}

\noindent This formula holds for any $u,v\in TM$. In our case 
$\{u,,v\}$ are orthonormal generators of $\Lambda$, so by (8) the third term on the right is  
$aJ(\alpha^{\#})$ where $J=J_{v\times u}$ is the complex structure of ${\bf V}$ given by Lemma 1 (and remarks following it). 

\vspace{.05in}

Now  if we call $w_{0}=\frac{1}{2}[ (u\times v) -w]$, and choose $\alpha $ among $1$-forms whose $ {\bf E}$ component zero (i.e. section of ${\bf V}$) with $|\alpha^{\#}|<1$ (hence $ a \neq 0 $),  the equation $(u\times v)_{\lambda}|_{Y}=w$ gives 
$w_{0}= |\alpha|^2 (u\times v) +a J(\alpha^{\#})$. By taking inner products of both sides with basis elements of ${\bf E}$, we see that $w_{0}^{\perp} = a J(\alpha^{\#}) $ where $w_{0}^{\perp} $ is the ${\bf V}$-component of $w_{0}$.  We can apply $J$ to both sides and solve $\alpha^{\#} = -\frac{1}{a} J( w_{0}^{\perp} )$. \qed 

\vspace{.1 in}

\subsection{Deforming associative submanifolds}

$\;$
\vspace{.1 in}

Let $G(3,7)\cong SO(7)/SO(3)\times SO(4)$ be the Grassmannian manifold of oriented $3$-planes in $\R^7$, and  $G^{\varphi_{0}}(3,7)=\{L \in G(3,7)\;| \;\varphi_{0}|_{L}=vol(L)\} $ be the submanifold of  associative $3$-planes. Recall that $G_2$ acts on  $G^{\varphi_{0}}(3,7)$ with stabilizer $SO(4)$ giving the identification
$G^{\varphi_{0}}(3,7)=G_2/SO(4)$ \cite{hl}. Recall also that if $\E \to G(3,7)$ and 
$\V =\E^{\perp}  \to G(3,4)$ are the canonical $3$-plane bundle and the complementary $4$-plane bundle, then we can identify the tangent bundle by 
$TG(3,7)=\E^{*}\otimes \V$. How does the tangent bundle of $G^{\varphi_{0}}(3,7)$ sit inside of this? The answer is given by the following Lemma.  By (7)  the cross product operation maps $\E \times \V\to \V$, and the metric gives an identification $\E^{*}\cong \E$, now if   $L=\langle e_1,e_2,e_3\rangle \in G^{\varphi_{0}}(3,7)$ with $\{ e_1,e_2, e_3=e_1\times e_2\}$ orthonormal, then

\vspace{.1 in}

{\Lem  $T_{L}G^{\varphi_{0}}(3,7)=\{\; \sum_{j=1}^{3} e^{j}\otimes v_{j} \in  \E^{*}\otimes \V \;|\; \sum e_j\times v_j =0  \;\}$.}

\proof  A tangent vector of  $G(3,7)$ at $L$ is a path of planes generated by  three orthonormal vectors $L(t)=\langle e_1(t),e_2(t),e_3(t)\rangle$, such that $L(0)=L$, in other words $\dot{L}=\sum e_{j}\otimes \dot{e}_{j}$. Clearly this tangent vector lies in $G^{\varphi_{0}}(3,7)$ if $e_{3}(t)=e_{1}(t)\times e_{2}(t)$. So $\dot{e}_{3}=\dot{e}_{1}\times e_{2}+ e_1 \times \dot{e}_{2}$. By taking cross product of both sides with $e_3$ and then using the identity (6) we get 
$$ \chi (\dot{e}_{1},e_2,e_3) + \chi (e_1, \dot{e}_{2}, e_{3}) +\chi (e_1,e_2,\dot{e}_{3})=0.$$
Now by using (8) and the fact that the cross product of two of the vectors in $\{e_1,e_2,e_3\}$ is equal to the third (in cyclic ordering), we get the result. \qed.

\vspace{.1 in}

It is easy to see that the normal bundle of $G^{\varphi_{0}}(3,7)$ in 
$G^(3,7)$ is isomorphic to $\V$ giving the exact sequence of the bundles over $G^{\varphi_{0}}(3,7)$:
$$ 0\to TG^{\varphi_{0}}(3,7) \to TG^(3,7) \stackrel {\times}\longrightarrow \V \to 0$$  

\vspace{.1 in}

From (7) we know that, if $Y^3\subset (M,\varphi )$ associative and $\nu$ is  its normal bundle, then the cross product operation maps:  $TY \times \nu \to \nu $, $\nu \times \nu \to TY$, and $TY \times TY \to TY $.
Let $\{ e_1,e_2,e_3 \}$ and $\{ e^1,e^2,e^3 \}$ be local frames and the dual coframes on $TY$ and ${\bf A_{0}}$ be the background Levi-Civita connection on $\nu$ induced from the metric on $M$ (there is also the identification $TY\cong T^{*}Y$ by  induced metric).  Then we can define a Dirac operator $\Di_{\bf A_{0}}:\Omega^{0}(\nu)\to \Omega^{0}(\nu)$ as the covariant derivative $\nabla_{\bf A_{0}}=\sum e^{j}\otimes \nabla_{e_{j}}$ followed by the cross product: 
 \begin{equation}
 \Di_{\bf A_{0}}= \sum e^{j}\times \nabla_{e_{j}} .
 \end{equation}

 
So the cross product plays the role of the Clifford multiplication in defining the Dirac operator in the normal bundle. We can extend this multiplication  to $2$-forms: $(a\wedge b) \times x=\frac{1}{2}[\;a\times (b\times x)-b\times (a\times x)\;]$  then by using (6) we get:
$$ (a \wedge b) \times x =\frac{1}{2} [\;x \lrcorner \;(a \wedge b) \;] -\chi(a,b,x). $$
In particular, when $a,b\in TY$ and $x\in \nu$ then   $(a \wedge b) \times x = -\chi(a,b,x) $. As usual we can twist this Dirac operator by connections on $\nu$, by replacing ${\bf A_{0}}$ with ${\bf A_{0}} +a$,  where $a\in \Omega^{1}(Y, ad \nu)$ is an endomorphism of $\nu$ valued $1$-form.
The following from \cite{as1}, is a generalized version of McLean's theorem \cite{m}. 

{\Thm The tangent space to associative 
submanifolds of a manifold with a $G_2$ structure $(M,\varphi )$ at an associative submanifold $Y$ is given by the kernel of the the twisted Dirac operator $\Di_{\bf A}: \Omega^{0}(\nu)\to \Omega^{0}(\nu)$, where ${\bf A}={\bf A_{0}}+a$ for some $a\in \Omega^{1}(Y, ad(\nu))$. The term $a=0$ when $\varphi$ is integrable.

\proof Recall the notations of Lemma 5. Let $L=\langle e_1,e_2,e_3 \rangle $ be  a tangent plane to $Y\subset M$. Any normal vector field $v$ to $Y$  moves $L$ by one parameter group of diffeomorphisms giving a path of $3$-planes in $M$, hence  it gives a vertical tangent vector 
$\dot{L}=\sum e_{j}\otimes \CL_{v}(e_j)\in T_{L}G_{3}(M)$ of the Grassmannian bundle of $3$-planes $G_{3}(M)\to M$ (where $\CL_{v}$ is the Lie derivative along $v$). By Lemma 5 this path of planes remain associative if  $\sum e_{j}\times \CL_{v}(e_j)=0$. Since  $\CL_{v}(e_j)=\bar{\nabla}_{e_j}(v)-\bar{\nabla}_{v}(e_j)$, where $\bar{\nabla}$ is the (torsion free) metric connection of $M$; then the result follows by letting $a(v)=\sum e^{j}\times \nabla_{v}(e_j)$ where $\nabla $ is the normal component of
$\bar{\nabla}$. If $\varphi $ is integrable, then  on a local chart it coincides with $\varphi_{0} $ up to quadratic terms, so $0=\nabla_{v}(\varphi)|_{Y}=\nabla_{v}(e^1\wedge e^2 \wedge e^3)$, which implies $a=0$. Also, by using the fact that the cross product operation preserves the tangent space of the associative manifold $Y$, it is easy to check that the expression of $a$ is independent of the choice of the orthonormal basis of $L$.\qed

\vspace{.05 in}

Notice that at any point by choosing normal coordinates we can make $a=0$. This reflects the fact that $\varphi $ coincides only pointwise with $\varphi_{0}$, not on a chart.  To make the Dirac operator onto, we can twist it by $1$-forms $a\in \Omega^{1}(Y)$, i.e. 

\vspace{.05 in}

{\Lem For associative $Y \subset (M,\varphi)$ the map $\Omega^{1}\times \Omega^{0}(\nu)\to  \Omega^{0}(\nu)$ defined by  $(x,a) \mapsto D_{\bf A }(x) + a\times x$ is onto, (by using appropriate Sobolev norms) }

\proof It suffices to show that the orthogonal complement of the image of this map is zero: Assume $\langle D_{\bf A }(x) , y\rangle +  \langle a\times x,y\rangle=0$ for all $x$ and $a$, then by taking $a=0$ and using the self adjointness of the Dirac operator we get   $D_{\bf A}(y)=0$. Hence $  \langle a\times x,y\rangle=0$, then the fact that the map $(x,a)\mapsto a\times x$ is surjective gives the result. Note that by (6) $a\times (a\times x)=- |a|^2 x$. \qed.

\vspace{.05 in} 
So for a generic choice of $a$ this twisted Dirac operator is onto, but what does this mean in terms of  the deformation space of the associative submanifolds? The  next Proposition (\cite{as1}) gives an answer. It says that if we perturb the deformation space with the gauge group (i.e. allowing a slight rotation of $TY$  by the gauge group of $TM$ during deformation) then it becomes smooth.  

\vspace{.05 in} 

Note that Theorem 6 may be explained by a  Gromov-Witten set-up:  Let  $G_{3}^{\varphi}(M)\subset G_{3}(M)\to M$ be the subbundle of  associative $3$-planes  with fiber $G_{3}^{\varphi}\cong G_2/SO(4)$ (\cite{hl}). We can form a bundle $G_{3}(Y, M)\to Im(Y^3,M)$ over the space of imbeddings, whose fiber over $f:Y\hookrightarrow M$ are the liftings $F$:
\begin{equation*}
\begin{array}{rrl}
  & & G_{3}(M) \supset G_{3}^{\varphi}(M)  \\
 F  \hspace{-.15in} &  \nearrow \;  &   \downarrow\\
\;\; Y\;\; & \stackrel{f}{\longrightarrow}   & M
\end{array}
\end{equation*}
The Gauss map $f\mapsto \sigma(f)$ gives a natural section to this bundle, and  $Y$ is associative if and only if this section maps into $G_{3}^{\varphi}(M)$. Theorem 6 gives the condition that the derivative of $\sigma$  maps into the tangent space of $G_{3}^{\varphi }(Y, M)$, which is the subbundle of $G_{3}(Y, M) $ consisting of $F$'s mapping into $G_{3}^{\varphi}(M)$. Recall that if $P\to M$ denotes the tangent frame bundle of $M$, then the gauge group $\CG(M)$ of $M$ is defined to be the sections of the $SO(7)$-bundle $Ad(P)\to M$, where $Ad(P)=P\times SO(7)/(p,h)\sim (pg, g^{-1}hg)$. By perturbing the Gauss map with the gauge group (i.e.  by composing $\sigma$  with the gauge group action $G_{3}(M)\to G_{3}(M)$ we can make it transversal to $G_{3}^{\varphi }(Y, M)$.

{\Prop The map $\tilde{\sigma}: \CG(M)\times Im(Y,M)\to G_{3}(Y,M)$ is transversal to  $G_{3}^{\varphi }(Y, M)$, where $\tilde{\sigma} (s,f)= s(f)\sigma(f)$ }

\proof We start with the local calculation of the proof of Theorem 6, except in this setting we need to take $\dot{L}=\sum e_{j}\otimes \CL_{v}(s e_{j})$, where $s\in SO(7)$ is the gauge group in the chart. Then the resulting equation is $\Di_{\bf A}(v) + \sum e_{j}\times {\bf v}(s) e_{j}=0$, where  ${\bf v}(s) e_{j}$ denotes the normal component of $v(s) e_{j}$ (here we are doing the calculation in normal coordinates where $\nabla_{v}(e_{k})=0$ pointwise). Then the argument as in the proof of Lemma 7 (by showing the second term is surjective) gives the proof.
\qed.

\vspace{.1 in}

  The kernel of this operator gives the deformations of  {\it pseudo associative submanifolds} of $Y\subset (M,\varphi)$, defined in  \cite{as1} as the manifolds where the perturbed Gauss map $\tilde{\sigma}$ maps $Y$ into $G_{3}^{\varphi}(M)$. We can either choose a generic $a$, or constraint the new variable $a$ by a natural second equation, which results equations resembeling the  Seiberg-Witten equations  as follows: 
  
\vspace{.05 in} 
  
  Let   $Y\subset (M,\varphi ,\Lambda)$ be a $\Lambda$-associative submanifold, we deform $Y$ in the  complex bundle $S\cong K^{-1}+\C = \W^{+}$ defined in Lemma 3.  From projections (9),  the  background $SO(4)$ metric connection on the normal bundle $\nu$  along with a choice of a connection $A$ on the line bundle $K\to Y$ gives a connection of the $Spin^{c}$ bundle, which in turn induces a connection  on the associated $U(2)$ bundle $W^{+}\to Y$. By using the Clifford multiplication  $TY_{\C}\otimes W^{+}\to W^{+}$ coming from the cross product (Lemma 3) we can form the  Dirac operator $\Di_{A}:\Omega^{0}(W^{+})\to \Omega^{0}(W^{+})$, whose kernel  identifies locally the deformations of $Y$ in the bundle $W^{+}$. Then if we constraint the  new variable $A$ by (11) we obtain deformations resembling  to the Seiberg-Witten  equations .
  \begin{equation}
\begin{array}{l}
  \Di_{\bf A}(x) =0 \\
*F_{A}=   \sigma(x).
\end{array}
  \end{equation}
\noindent where $*$ is the star operator of $Y$ induced from the background metric of $M$, 
 and $(x,A)\in \Omega^{0}(W^{+})\times \CA(K)$, and $\CA(K)$ is the space of connections on $K$.  From Weitzenb\"{o}ck formula and (12), the above equations give compactness to this type of local deformation space, hence allow us to assign Seiberg-Witten invariant to $Y$. Now the natural question is how easy to produce $\Lambda$-associative submanifolds $Y^{3}\subset (M,\varphi ,\Lambda)$? One answer is that any zero set $Y^3$ of a transverse section of  ${\bf V}\to M$ gives a $\Lambda$-associative submanifold. This is because the transversality gives a canonical identification $TY\cong {\bf E}|_{Y}$. Then the natural question is: Are there natural sections of ${\bf V}$? We can obtain such things from the other $G_2$ structures as follows. Recall that $\Lambda=\langle u,v\rangle$ gives the section $s= \langle u,v,u\times v\rangle$ of the bundle $G_{3}(M)\to M$. For any other $G_2$ structure $\psi$  on $(M,\varphi)$ defines a section $\chi_{\psi}$ of $\V \to G_{3}(M)$. Then by pulling $\chi_{\psi}$ with $s$ over $M$ produces a natural section ${\bf s} (\varphi,\psi)$ of ${\bf V}\to M$. To these sections we can associate an integer valued invariant, i.e. the Seiberg-Witten invariant of their zero set (15) (the nontransverse sections we can associate zero). Consequences of this will be explored in a future paper. 

\vspace{.05in}
 
 Note that in the usual Seiberg-Witten equations on  Y, we use an action of $T^{*}Y$: $W^{+}\to W^{-}$ coming from $Spin^{c}$ structure, which then extends an action of $\Lambda^{2}(Y)$ to  $W^{+}\to W^{+}$. Here $T^{*}Y$  acts as $W^{+}\to W^{+}$ by Lemma 3. On the other hand by the background metric we have the identification $T^{*}Y\cong \Lambda^{2}(Y)$.

\end{document}